\documentclass[12pt]{article}
\input{xypic}
\usepackage{amsmath}
\usepackage{amssymb}
\usepackage{latexsym}
\usepackage{enumerate}

\newtheorem{theorem}{Theorem}[section]
\newtheorem{lemma}[theorem]{Lemma}

\newtheorem{proposition}[theorem]{Proposition}
\newtheorem{corollary}[theorem]{Corollary}

\newenvironment{proof}{\noindent{\sc Proof.}}{\hfill$\square$}

\numberwithin{equation}{section}

\newcommand{\A}{{\mathcal A}}
\newcommand{\mO}{{\mathcal O}}
\newcommand{\U}{{\mathcal U}}

\newcommand{\C}{\mbox{$\mO$-alg}}

\newcommand{\Q}{{\mathbb Q}}
\newcommand{\epi}{\mbox{$\to$\hspace{-0.35cm}$\to$}}
\newcommand{\mono}{\hookrightarrow}

\begin{document}

\title{Plus-construction of algebras over an operad,
cyclic and Hochschild homologies up to homotopy
\thanks{Primary: 19D06, 19D55; Secondary: 18D50, 18G55, 55P60, 55U35}}

\author{
{\sc David Chataur, Jos\'e L. Rodr\'{\i}guez, and J\'er\^ome
Scherer}
\thanks{The first author was supported by Marie Curie grant HPMF-CT-2001-01179,
the second by EC grant HPRN-CT-1999-00119, CEC-JA grant FQM-213,
and DGIMCYT grant BFM2001-2031.}}

\date{\today}
\maketitle

\begin{abstract}
In this paper we generalize the plus-construction given by M.
Livernet for algebras over rational differential graded Koszul
operads to the framework of admissible operads (the category of
algebras over such operads admits a closed model category
structure). We follow the modern approach of J. Berrick and C.
Casacuberta defining topological plus-construction as a
nullification with respect to a universal acyclic space.
Similarly, we construct a universal $H^Q_*$-acyclic algebra
$\mathcal U$ and we define $A\longrightarrow A^+$ as the $\mathcal
U$-nullification of the algebra $A$. This map induces an
isomorphism on Quillen homology and quotients out the maximal
perfect ideal of $\pi_0(A)$. As an application, we consider for
any associative algebra $R$ the plus-constructions of $gl(R)$ in
the categories of Lie and Leibniz algebras up to homotopy. This
gives rise to two new homology theories for associative algebras,
namely cyclic and Hochschild homologies up to homotopy. In
particular, these theories coincide with the classical cyclic and
Hochschild homologies over the rational.
\end{abstract}

\section*{Introduction}

Quillen's plus construction for spaces was designed so as to yield
a definition of higher algebraic $K$-theory groups of rings.
Indeed, for any $i \geq 1$, $K_i R = \pi_i BGL(R)^+$, where
$GL(R)$ is the infinite general linear group on the ring $R$. The
study of the additive analogue, namely the Lie or Leibniz algebra
$gl(R)$ has already produced a number of papers showing the strong
link with cyclic and Hochschild homology (For classical background
on these theories we refer to \cite{LOC} and to the survey
\cite{LOS}). However there have always been restrictions, such as
working over the rationals.

For example M. Livernet has given a plus-construction for algebras
over a Koszul operad in the rational context \cite{ML} by way of
cellular techniques imitating the original topological
construction given by D. Quillen (see also \cite{Pi} for a
plus-construction in the context of simplicial algebras).
Specializing to the category of Lie, respectively Leibniz
algebras, she proved then that the homotopy groups of $gl(R)^+$
are isomorphic to the cyclic, respectively Hochschild homology
groups of $R$.

In the category of topological spaces plus-construction can be
viewed as a localization functor, which has the main advantage to
be functorial. This idea goes probably back to A.K.~Bousfield and
E.~Dror Farjoun, but the first concrete model was given by
J.~Berrick and C.~Casacuberta in \cite{BC}. They provide a
``small" universal acyclic space $B\mathcal F$ such that the
nullification $P_{B\mathcal F} X$ is the plus-construction $X^+$.

Recently, thanks to the work of P.~Hirschorn \cite{Hir} it appears
possible to do homotopical localization in a very general
framework. In fact, one can construct localizations in any closed
model category satisfying some mild extra conditions (left proper
and cofibrantly generated), such as categories of algebras over
admissible operads. The category of Lie algebras over an arbitrary
ring is not good enough for example. One needs to take first a
cofibrant replacement $\mathcal L_\infty$ of the Lie operad and
can perform localization in the category of $\mathcal
L_\infty$-algebras, which we call Lie algebras up to homotopy.

This allows to define a functorial plus-construction in the
category of algebras over an admissible operad as a certain
nullification functor with respect to an algebraic analogue $\U$
of Berrick's and Casacuberta's acyclic space. This extends the
results of M. Livernet to the non-rational case. In the following
theorem $\mathcal P\pi_0(A)$ denotes the maximal $\pi_0
\mO$-perfect ideal of $\pi_0 A$.

\medskip

\noindent {\bf Theorem \ref{plus}}.
{\it Let $\mathcal O$ be an admissible operad. Then the
homotopical nullification with respect to $\mathcal U$ is a
functorial plus-construction in the category of algebras over
$\mathcal O$. It enjoys the following properties:}
\begin{itemize}
\item[(i)]{} $A^+ \simeq Cof(\coprod_{f\in Hom(\mathcal U,A)} \U
\longrightarrow A)$

\item[(ii)]{} $\pi_0(A^+)\cong\pi_0(A)/\mathcal P\pi_0(A)$

\item[(iii)]{} $H^Q_*(A)\cong H^Q_*(A^+)$
\end{itemize}

\medskip

Of particular interest is the plus-construction in the category of
Lie algebras up to homotopy. If we apply these constructions to
the algebra $gl(R)$ of matrices of an associative algebra, and if
we consider it as a Lie algebra up to homotopy, we obtain what we
call the cyclic homology theory up to homotopy. Thus $HC_i^\infty
(R)$ is defined as $\pi_i gl(R)^+$ for any $i \geq 0$. This theory
corresponds to the classical cyclic homology over the rationals.
We summarize in a proposition the computations of the lower
homology groups (see Proposition~\ref{K0}, \ref{K1}, \ref{K2}).
They share a striking resemblance with the low dimensional
algebraic $K$-groups ($st(R)$ stands for the Steinberg Lie
algebra, see Section~5).

\medskip

\noindent {\bf Proposition}.
{\it Let $k$ be a field and $R$ be an associative $k$-algebra.
Then
\begin{itemize}
\item[(1)] $HC_0^{\infty}(R)$ is isomorphic to $R/[R,R]$.

\item[(2)] $HC_1^{\infty}(R)$ is isomorphic to $Z(st(R)) \cong
H_1^Q(sl(R))$.

\item[(3)] $HC_2^{\infty}(R)$ is isomorphic to $H_2^Q(st(R))$.

\item[(4)] $HC_3^{\infty}(R)$ is isomorphic to $H_3^Q(st(R))$.
\end{itemize}}

The same kind of results hold for Hochschild homology up to
homotopy, which is defined similarly using Leibniz algebras. The
above computation in the case $k = \Q$ yields the following
corollary for classical cyclic homology.

\medskip

\noindent {\bf Corollary \ref{classical}}.
{\it Let $R$ be an associative $\mathbb Q$-algebra. Then
\begin{itemize}
\item[(1)] $HC_0(R)$ is isomorphic to $R/[R,R]$.

\item[(2)] $HC_1(R)$ is isomorphic to $Z(st(R)) \cong
H_2^{Lie}(sl(R))$.

\item[(3)] $HC_2(R)$ is isomorphic to $H_3^{Lie}(st(R))$.

\item[(4)] $HC_3(R)$ is isomorphic to $H_4^{Lie}(st(R))$.
\end{itemize}}

Computation $(1)$ is well known and trivial, whereas $(2)$, $(3)$
and $(4)$ are non trivial results (for $(2)$ and $(3)$ we refer to
\cite {KL}). The plan of the paper is as follows. First we
introduce the notion of algebra over an operad and recall when and
how one can do homotopy theory with these objects. In a short
second section we explain what homotopical localization and
nullification functors are for algebras. The main theorem about
the plus-construction appears then in Section~3 and Section~4
contains the properties of plus-construction with respect to
fibrations and extensions. The final section is devoted to the
computations of the low dimensional additive $K$-theory groups.

\medskip

\noindent {\bf Acknowledgements.} We would like to thank Benoit
Fresse and Jean-Louis Loday for helpful comments. The first author
thanks Muriel Livernet for sharing so many ``operadic" problems,
her thesis was the starting of all this story. The first two
authors would like to thank Yves F\'elix and the Universit\'e
Catholique de Louvain La Neuve for having made this collaboration
possible. We also thank the Centre de Recerca Matematica, the
University of Almer\'\i a, the Universit\'e de Lausanne, as well
as the Universitat Autonoma de Barcelona for their hospitality.

\section{Operads and algebras over an operad}
We fix $R$ a commutative and unitary ring. We work in the category
$\bf R-dgm$ of differential $\mathbb Z$-graded $R$-modules and
especially with chains (the differential decreases the degree by
$1$). All the objects we will consider in this article will be in
fact $\mathbb N$-graded, but for technical reasons, namely in
order to use techniques of fiberwise localization, it is handy to
view them as unbounded chain complexes. For classical background
about operads and algebras over an operad see \cite {GJ}, \cite
{GK}, \cite{KM} and \cite{Lo}.
\\
\\
{\bf $\Sigma_*$-modules.} A $\Sigma_*$-module is a sequence
$\mathcal M=\{\mathcal M(n)\}_{n>0}$ of objects $\mathcal M(n)$ in
$\bf R-dgm$ together with an action of the symmetric group
$\Sigma_n$. The category of $\Sigma_*$-modules is a monoidal
category. We denote by $\mathcal M\circ \mathcal N$ the product of
two $\Sigma_*$-modules and by $\mathbf 1$ the unit of this
product. The unit is defined by $\mathbf 1(1)=R$ and $\mathbf
1(i)=0$ for $i \not = 1$.
\\
\\
{\bf Operads.} An {\bf operad} $\mathcal O$ is a monad in the
category of $\Sigma_*$-modules. Hence we have a product
$\gamma:\mathcal O\circ \mathcal O \longrightarrow \mathcal O$
which is associative and unital. Equivalently the product $\gamma$
defines a family of composition products
$$
\gamma:\mathcal O(n)\otimes\mathcal
O(i_1)\otimes\ldots\otimes \mathcal O(i_n)\longrightarrow \mathcal
O(i_1+\ldots+i_n)
$$
which must satisfy equivariance, associativity and unitality
relations (also called May's axioms). Moreover we suppose that
$\mathcal O(1)=R$ and that each chain complex $\mathcal O(n)$ is
concentrated in positive differential degrees (i.e. $\mathcal
O(n)_p=0$ for any $p<0$).
\\
There is a free operad functor:
$$F:\Sigma_*-modules\longrightarrow Operads$$
which is left adjoint to the forgetful functor. It can be defined
using the formalism of trees.
\\
\\
{\bf Algebras over an operad.} Let us fix an operad $\mathcal O$.
An algebra over $\mathcal O$ (also called $\mathcal O$-algebra) is
an object $A$ of $\bf R-dgm$ together with a collection of
morphisms
$$
\theta:\mathcal O(n)\otimes A^{\otimes
n}\longrightarrow A
$$
called evaluation products, which are equivariant, associative and
unital.
\\
There is a free $\mathcal O$-algebra functor
$$S(\mathcal O,-):\bf R-dgm\longrightarrow \mathcal O-algebras$$
which is left adjoint to the forgetful functor. For any $M\in \bf R-dgm$
it is given by $S(\mathcal O,M)=\bigoplus_{n>0} \mathcal
O(n)\otimes_{R[\Sigma_n]}M^{\otimes n}$.
\\
\\
{\bf Classical operads.}
\\
a) Let $M$ be an object of $\bf R-dgm$. We associate to it the
endomorphism operad given by:
$$
End(M)(n)=Hom_{\bf R-dgm}(M^{\otimes n},M).
$$
Any $\mathcal O$-algebra structure on $M$ is given by a morphism
of operads $\mathcal O\longrightarrow End(M).$
\\
b) The operad $\mathcal Com$, defined by $\mathcal Com(n)=R$. The
$\mathcal Com$-algebras are the differential graded algebras.
\\
c) The operad $\mathcal As$ defined by $\mathcal
As(n)=R[\Sigma_n]$. The $\mathcal As$-algebras are precisely the
differential graded associative algebras.
\\
d) The operad $\mathcal Lie$. A $\mathcal Lie$-algebra $L$ is an
object of $\bf R-dgm$ together with a bracket which is
anticommutative and satisfies the Jacobi relation. If $2\in R$ is
invertible then a $\mathcal Lie$-algebra is a classical Lie
algebra. Otherwise, the category of classical Lie algebras appears
as full subcategory of the category of $\mathcal Lie$-algebras.
\\
e) The operad $\mathcal Leib$ which is the operad of Leibniz
algebras. Leibniz algebras are algebras equipped with a bracket of
degree zero that satisfies the Jacobi relation. A Lie algebra is
an anti-commutative Leibniz algebra. We have an epimorphism of
operads
$$\mathcal Leib\longrightarrow\mathcal Lie.$$
\\
{\bf Homotopy of operads.} In \cite {Hi1}, and \cite {BM},  V.
Hinich and C. Berger-I. Moerdijk proved that the category of
operads is a closed model category. This structure is obtained
from the one on the category $\bf R-dgm$ via the free operad
functor. In this category the weak equivalences are the
quasi-isomorphisms and the fibrations are the epimorphisms. The
cofibrant operads are the retracts of the quasi-free operads.
\\
\\
{\bf Homotopy of algebras over an operad.}
Let $W_d$ be the following object of $\bf R-dgm$:
$$\ldots 0\longrightarrow R=R \longrightarrow 0\ldots$$
concentrated in differential degrees $d$ and $d+1$ ($d\in \mathbb
Z$). Using the terminology of \cite {BM}, we say that $\mathcal O$
is admissible if the canonical morphism of $\mathcal O$-algebras:
$$
A\longrightarrow A\coprod S(\mathcal O,W_d)
$$
is a quasi-isomorphism for any $\mathcal O$-algebra $A$ and for
all $d$. For any admissible operad $\mathcal O$ there exists a
closed model structure on the category of $\mathcal O$-algebras,
which is transferred from {\bf R-dgm} along the free-forgetful
adjunction given by $S(\mO, -)$.
As for operads the weak equivalences are the
quasi-isomorphisms, the fibrations are the epimorphisms, and the
cofibrant $\mathcal O$-algebras are the retracts of the quasi-free
$\mathcal O$-algebras.
\\
The category of $\mathcal O$-algebras is cofibrantly generated and
cellular in the sense of P. Hirsch\-horn \cite {Hir}. The set of
generating cofibrations is
$$
I=\{i_n:\mathcal O(x_n)\longrightarrow \mathcal O(x_n,y_{n+1})\}
$$
where $\mathcal O(x_n)$ is the free $\mO$-algebra on a generator of
degree $n$ and $\mathcal O(x_n,y_{n+1})$ is the free $\mO$-algebra over
the differential graded module $R<x_n,y_{n+1}>$ with two copies of
$R$, one in degree $n$ the other in degree $n+1$, the differential
of $y_{n+1}$ being $x_n$. The set of generating acyclic
cofibrations is
$$
J=\{j_n:0\longrightarrow \mathcal O(x_n)\}.
$$
Notice that the free algebra $\mathcal O(x_n)$ plays the role of
the sphere $S^n$.
\\
\\
Over the rational numbers all operads are admissible. This is not
the case over an arbitrary ring, for example the operads $\mathcal
Com$ and $\mathcal Lie$ over the integers are not. However
cofibrant operads and the operad $\mathcal As$ are always
admissible.
\\
\\
{\bf Homology of algebras.} Let $\mathcal O$ be an admissible
operad, and let $A$ be an $\mathcal O$-algebra. An element $a\in
A$ is called decomposable if it lies in the ideal $A^2$, the image
of the evaluation products
$$
\theta(n):\mathcal O(n)\otimes A^{\otimes n}\longrightarrow A
 $$
for any $n>1$. We denote by $QA=A/A^2$ the space of
indecomposables of the algebra~$A$. The Quillen homology of $A$,
denoted by $H^Q_*(A)$, is the homology of $QS(\mathcal O,V)$ where
$S(\mathcal O,V)$ is a cofibrant replacement of $A$. This does not
depend on the choice of the cofibrant replacement.
\\
Moreover, any cofibration sequence $A\longrightarrow
B\longrightarrow C$ of $\mathcal O$-algebras yields a long exact
sequence in Quillen homology.
\\
In \cite{BF} Fresse generalizes Koszul duality for operads as
defined by Ginzburg and Kapranov in \cite{GK}. As the operads
$\mathcal Lie$ and $\mathcal Leib$ are Koszul operads, over
$\mathbb Q$ one can compute their Quillen homology by way of a
nice complex:
\\
{\bf $\mathcal Lie$-algebras.} Let $L$ be a $\mathcal Lie$-algebra
over $\mathbb Q$. The homology of $L$, denoted here by
$H^{Lie}_*(L)$, is computed using the Chevalley-Eilenberg complex
$CE_*(L)$. Now we can consider $L$ as a $\mathcal
L_{\infty}$-algebra and compute $H^Q_*(L)$. Fresse's results give
the following isomorphism:
$$H^Q_*(L)\cong H^{Lie}_{*+1}(L).$$

\noindent {\bf $\mathcal Leib$-algebras.} The same kind of results
hold for $\mathcal Leib$-algebras.  Consider a Leibniz algebra
$L$. The homology of $L$, denoted by $H^{Leib}_*(L)$, is computed
using a complex described in \cite{Li2}. One has again a similar
isomorphism:
$$H^Q_*(L)\cong H^{Leib}_{*+1}(L).$$

\noindent {\bf Quillen cohomology of discrete algebras.} We refer
the reader to \cite {GH} for more details and precise definitions.
A discrete algebra is an $\mathcal O$-algebra concentrated in
differential degree $0$. The structure of $\mathcal O$-algebra
reduces then in fact to a structure of $\pi_0(\mathcal
O)$-algebra.
\\
In the case of discrete algebras there is also a notion of Quillen
cohomology with coefficents. Fix a discrete $\mathcal O$-algebra
$A$ and a discrete $A$-module $M$.
\\
A derivation $D: A \rightarrow M$ in $Der(A,M)$ is a linear map
such that for any $o\in\mathcal O(n)$ we have:
$$
D(o(a_1,\ldots,a_n))=\sum_{i=1}^n o(a_1,\ldots,D(a_i),\ldots,a_n).
$$
We can define Quillen cohomology by computing the derived functors
of $Der(A,M)$, that is by taking $A'$ a cofibrant replacement of
$A$ in the category of $\mathcal O$-algebras and computing the
homology of the complex $Der(A',M)$. This has also a homotopical
interpretation: Consider $M$ as a trivial $\mathcal O$-algebra and
denote by $\Sigma^n M$ the $n$-th suspension of $M$ in the
category ${\bf R-dgm}$. Then
$$
H^n_Q(A,M)=[A,\Sigma^n M]_{\mathcal O-alg}\cong [QA,\Sigma^n
M]_{\bf R-dgm}. $$
Moreover $H^1_Q(A,M)$ classifies square zero
extensions of $A$ by $M$. A square zero extension is an exact
sequence of modules:
$$
0\rightarrow M\rightarrow B\stackrel{p}{\rightarrow} A\rightarrow 0,
$$
such that $p$ is a morphism of $\mathcal O$-algebras and the
structure of $\mathcal O$-algebra on $M$ is trivial. The set of
isomorphism classes of square zero extension is denoted by
$Ex(A,M)$. By a classical result of Quillen \cite {QU} we have the
following isomorphism:
$$
H^1_Q(A,M)\cong Ex(A,M).
$$

\noindent{\bf Hurewicz Theorem.} In her thesis
\cite[Theorem~2.13]{ML} M. Livernet proved a Hurewicz type theorem
for algebras over a Koszul operad in the rational case. A result
of Getzler and Jones about the construction of a cofibrant
replacement for $\mathcal O$-algebras, which uses the Bar-Cobar
construction, extends the proof of Livernet to admissible operads.

\begin{theorem}
\label{Hurewicz}
Let $A$ be an $\mathcal O$-algebra. We suppose that the underlying
chain complex of $A$ is concentrated in non-negative degrees. Then
there is a Hurewicz morphism:
$$
Hu:\pi_*(A)\longrightarrow H_*^Q(A)
$$
induced by the projection on indecomposable elements. It satisfies
the following properties:
\\
i) If $\pi_k(A)=0$ for $0\leq k \leq n$ then $Hu$ is an
isomorphism for $k\leq 2n+1$ and an epimorphism for $k=2n+2$.
\\
ii) If $\pi_0(A)=0$ and $H_k^Q(A)=0$ for $0\leqq k \leqq n$ then
$Hu$ is an isomorphism for $k\leq 2n+1$ and an epimorphism for
$k=2n+2$.
\end{theorem}

\begin{proof}
In the case of $0$-connected chain complexes we have a Quillen
adjunction between $\mathcal O$-algebras and $B\mathcal
O$-coalgebras, \cite {GJ} and \cite {AC}. These two functors
provide for any $A$ a cofibrant replacement of the form
$S(\mathcal O,C(B\mathcal O,A))$ where $C(B\mathcal O,A)$ is the
coalgebra over the cooperad $B\mathcal O$ obtained by applying the
operadic bar construction. Now Livernet's arguments apply to
$C(B\mathcal O,A)$.
\end{proof}
\\
\\
\noindent {\bf Perfect algebras over an operad.} Consider an
algebra $A$ over an operad in the category of $R$-modules. The
algebra $A$ is called $\mathcal O$-perfect if any element in $A$
is decomposable i.e. $A=A^2$ or $QA=0$. We define ${\mathcal P}A$,
the maximal $\mathcal O$-perfect ideal of $A$, by transfinite
induction.
 \\
Let $A_0$ be the ideal $A^2$. We define the ideals $A_{\alpha}$
inductively by setting $A_{\alpha}=(A_{\alpha-1})^2$ if $\alpha$
is a successor ordinal and $\mathcal A_{\alpha}=\cap_{\beta <
\alpha}A_{\beta}$ if $\alpha$ is a limit ordinal. Then we set
$\mathcal PA=lim_{\alpha} A_{\alpha}$. This inverse system
actually stabilizes for some ordinal $\beta$, hence
$A_{\beta}^2=A_{\beta}$ and $\mathcal P A=A_{\beta}$. Of course,
if $QA=0$ then we have $\mathcal P A=A$.
 \\
Consider an epimorphism $f:\mathcal O\longrightarrow \mathcal O'$
of operads and let $A$ be an $\mathcal O'$-algebra. Then if $A$ is
$\mathcal O'$-perfect it is also $\mathcal O$-perfect. Thus we
have an inclusion $\mathcal P A\subseteq\mathcal P' A$ of the
$\mO$-perfect ideal into the $\mO'$-perfect ideal of $A$.
 \\
For any differential graded algebra $A$ over an operad $\mathcal
O$, denote the $i$-th homology group of the underlying chain
complex by $\pi_i(A)$. Then $\pi_0(A)$ is a $\pi_0(\mathcal
O)$-algebra in the category of $R$-modules. Moreover we have
$H^Q_0(A)=Q\pi_0(A)$.

\section{Homotopical localization and nullification}
The theory of homological and homotopical localization of
topological spaces developed by Bousfield and Dror Farjoun (see
e.g. \cite{Bou97}, \cite{Dro}) has an analogue in the category of
all algebras over an admissible operad $\mO$. This takes place in
the more general framework established by Hirschhorn in
\cite{Hir}. Our category $\C$ of all algebras over an admissible
operad $\mO$ is indeed cellular and left proper. We explain first
how to build mapping spaces in a model category which is not
supposed to be simplicial, and recall then what is meant by
homotopical localization with respect to a morphism in this
context.
\\
\\
{\bf Mapping spaces.} One way to construct mapping spaces up to
homotopy in a model category is to find a cosimplicial resolution
of the source (as in \cite[Definition 18.1.1]{Hir}). An elementary
method to get one is provided by \cite{CS}. Let $X$ be a cofibrant
$\mathcal O$-algebra and define $\Delta[n] \otimes X$ as the
homotopy colimit of the constant diagram with value $X$ over the
simplex category $\Delta[n]$. This is simply the category of
non-empty subsets of $\underbar n = \{1, \dots , n \}$ with
terminal object $\underbar n$. For example $\Delta[1] \otimes X$
is the homotopy colimit of the constant diagram $X \rightarrow X
\leftarrow X$. A cofibrant replacement of this diagram is given by
$X \rightarrow Cyl(X) \leftarrow X$ where $Cyl(X)$ is a cylinder
object for $X$ (by definition the folding map factorizes as a
cofibration followed by a weak equivalence $X \coprod X
\hookrightarrow Cyl(X) \stackrel{\simeq}{\twoheadrightarrow} X$).
Having a functorial cofibrant replacement guarantees that
$\Delta[-] \otimes X$ is a cosimplicial object. Therefore one can
define $map(X, Y) = mor_{\mathcal O-alg}(\Delta[-] \otimes X, Y)$.

Let $\partial \Delta[n] \otimes X$ be the homotopy colimit of the
constant diagram with value $X$ over the simplex category of the
boundary $\partial \Delta[n]$, i.e. the category $\Delta[n]$
without the terminal object. Define $S^n \wedge X$ to be the
homotopy cofiber of $\partial \Delta[n] \otimes X \hookrightarrow
\Delta[n] \otimes X$. For example $S^1 \wedge X$ is simply the
homotopy cofiber of $X \coprod X \hookrightarrow Cyl(X)$.

\begin{lemma}
For any cofibrant $\mathcal O$-algebra $X$ we have $S^1 \wedge X
\simeq \Sigma X$.
\end{lemma}

\begin{proof}
It is always true in a pointed model category that the homotopy
cofiber of $X \coprod X \hookrightarrow Cyl(X)$ is equivalent to
the suspension of $X$.
\end{proof}

\begin{proposition}
\label{smash}
In the category of $\mathcal O$-algebras $S^n \wedge X \simeq
\Sigma^n X$.
\end{proposition}

\begin{proof}
The proof is by induction on $n$, the case $n=1$ having been
proved in the preceding lemma. Recall that $\partial \Delta[n]
\otimes X$ can be computed by decomposing the cubical homotopy
colimit as a push-out involving homotopy colimits on the front and
back face (see for example \cite[Lemma~0.2]{Goo}): $\partial
\Delta[n] \otimes X \simeq hocolim( X \leftarrow \partial
\Delta[n-1] \otimes X \rightarrow X )$. Commuting homotopy
colimits again we see that $S^n \wedge X$ is the homotopy push-out
of the diagram $(* \leftarrow S^{n-1} \wedge X \rightarrow * )$,
i.e. $\Sigma ( S^{n-1} \wedge X ) \simeq \Sigma^n X$.
\end{proof}
\\
\\
{\bf Homotopy groups of mapping spaces.} When $X$ is cofibrant and
$Y$ fibrant, the simplicial set $map(X, Y)$ is fibrant, so its
homotopy groups can be computed. There is always at least a
morphism $X \rightarrow Y$, namely the trivial one. Let us denote
by $map^0 (X, Y)$ the component of this trivial map.

\begin{proposition}
\label{homotopygroups}
Let $X$ be a cofibrant and $Y$ a fibrant $\mathcal O$-algebra.
Then $\pi_n map^0 (X, Y) \cong [\Sigma^n X, Y]$.
\end{proposition}

\begin{proof}
The mapping space $map(X, Y)$ is fibrant, so that the $n$-th
homotopy group consists of homotopy classes of those $n$-simplices
whose faces are trivial. An $n$-simplex is a morphism $f:
\Delta[n] \otimes X \rightarrow Y$ and its faces are trivial if
the composite $\partial\Delta[n] \otimes X \hookrightarrow
\Delta[n] \otimes X \rightarrow Y$ is so. Therefore $f$ factorizes
through $S^n \wedge X$ and we conclude by Proposition~\ref{smash}.
\end{proof}
\\
\\
{\bf $f$-localization.} Let $\mO$ be an admissible operad. Let $f:
X\to Y$ be a morphism between two cofibrant $\mO$-algebras. An
$\mO$-algebra $Z$ is called {\it $f$-local} if the map of
simplicial sets
$$
map(f,Z): map(Y,Z)\to map(X,Z)
$$
is a weak homotopy equivalence. A morphism of $\mO$-algebras
$h:A\to B$ is called an {\it $f$-equivalence} if it induces a weak
homotopy equivalence $map(h,Z): map(B,Z)\to map(A,Z)$ for every
$f$-local algebra $Z$. Theorem $4.1.1$ from \cite{Hir} ensures
then the existence of an $f$-localization functor, i.e. a
continuous functor $L_f: \C\to \C$ together with a natural
transformation $\eta: Id\to L_f$ from the identity functor to
$L_f$, such that $\eta_A: A\to L_f A$ is an $f$-equivalence and
$L_f A$ is $f$-local for any $\mO$-algebra $A$.
\\
\\
{\bf $X$-nullification.} In the special case when $f$ is of the
form $f: X\to 0$ or $f:0\to X$, the functor $L_f$ is also denoted
by $P_X$, and called $X$-nullification functor. Note that an
algebra $Z$ is $f$-local or $X$-null if $map(X,Z)$ is weakly
homotopy equivalent to a point. By
Proposition~\ref{homotopygroups} this is equivalent to requiring
that $[\Sigma^k X, Z]$ be trivial for all $k\geq 0$.
\\
\\
{\bf Example:} An interesting example is when $X$ is the free
$\mO$-algebra $\mO(x)$ with one generator $x$ in dimension $n$.
This plays the role of the $n$-dimensional sphere, hence
$\mO(x)$-nullification gives rise to a functorial $n$-Postnikov
section in this category.

\section{An additive plus-construction}

A {\it Quillen plus-construction} of an algebra $A$ over an operad
$\mO$ is a Quillen homology equivalence $\eta: A\to A^+$ which
quotients out the perfect radical on $\pi_0$, that is
$$
\pi_0(A^+)\cong \pi_0(A)/ \mathcal P \pi_0(A).
$$
If it exists then it is unique up to quasi-isomorphism. This is
because of the following universal property: For any morphism
$g:A\to B$ with $Q\pi_0(B)=\pi_0(B)$, there exists a map $\tilde
g: \A^+ \to B$ such that $\tilde g \eta = g$, and which is unique
up to homotopy.

We now construct a functorial Quillen plus-construction as a
nullification with respect to a universal acyclic algebra, i.e. an
acyclic algebra $\U$ such that the associated nullification $A\to
P_{\U} A$ is the plus-construction.

\medskip

\noindent {\bf A universal acyclic algebra.} The algebra $\U$ is
defined as a big coproduct
$$
\U=\coprod_ {(T,\phi)} \U_{(T,\phi)},
$$
where $(T,\phi)$ ranges over all $\mO$-trees. Every
$\U_{(T,\phi)}$ will be the homotopy colimit of a certain direct
system of free $\mO$-algebras $\{U_r, \varphi_r\}_{r\geq 0}$
associated to $(T,\phi)$ in a canonical way. A {\it rooted tree}
$T=\{V(T),A(T)\}$ is a directed graph such that any vertex $v\in
V(T)$ has one ingoing arrow $a_v\in A(T)$, except the root that
has no ingoing arrow, and such that the following additional
conditions are satisfied: Each vertex $v$ has a finite number of
outgoing arrows, denoted by $val(v)$; the set $suc(v)$ of
successors vertices of $v$, i.e. those which are connected to $v$
by an ingoing arrow is finite and totally ordered; and finally,
the vertices $v$ of even level have at least 2 successors. The
root has level 0, and inductively we say that a vertex $v$ has
level $k$ if $v\in suc(u)$ for some $u$ of level $k-1$.

Let $\mO$ be any operad. An $\mO$-tree is a pair $(T,\phi)$ where
$T$ is a rooted tree and $\phi$ is a function which associates to
each vertex $v$ of odd level a multilinear operation $o_n\in
\mO(n)_0$ where $n=val(v)\geq 2$.

We next define the direct system $\{U_r, \phi_r\}$ of free
$\mO$-algebras associated to a given $\mO$-tree $(T,\phi)$, by
induction on $r$: Let $U_0$ be the free $\mathcal O$-algebra on
one generator $x$ in dimension 0 (corresponding to the root). Let
$n=val(root)$ and $suc(root)=\{v_1, \ldots, v_n\}$. For each $j=1,
\ldots, n$, let $k_j= val(v_j)$ and $o_{k_j}=\phi(v_j)$ be the
multilinear operation in $\mO(k_j)_0$ associated to the vertex
$v_j$. Choose $k_j$ free generators $x_{1j1}, x_{1j2}, \cdots,
x_{1jk_j}$ in dimension $0$ (corresponding to the vertices in
$suc(v_j)$ of level $2$, the first index indicates half of the
level). Let $U_1$ be the free $\mathcal O$-algebra on those
generators. Define $\phi_1: U_0\to U_1$ on the generator $x$ by
$$
\phi_1(x)= \sum_{j=1}^n \theta(o_{k_j};x_{1j1},x_{1j2},\ldots
,x_{1jk_j}).
 $$
Inductively, we define then $U_r$ as the free $\mO$-algebra on as
many generators as there are vertices of level $2r$, and $\phi_r:
U_{r-1}\to U_r$ is given on each generator of $U_{r-1}$ by a
similar formula as the above one for $\phi_1(x)$. Hence the
homotopy colimit $\U_{(T,\phi)}$ is free on generators $x_{I}$ of
degree $0$ and $y_{I}$ of degree $1$ where $I$ is a multi-index of
the form $ljs$, $l$ indicating half of the level where these
generators are created, $1 \leq s \leq k_j$, and the differential
$d(y_{I}) = x_{I} - \phi_l(x_{I})$.

\medskip

\noindent {\bf The cone of $\U$.} In order to do some computations
with this acyclic algebra, we need to describe how to construct
the cone of it. Let us simply describe the cone on $\U_{(T,\phi)}$
for a fixed tree~$T$. For each generator $x_I$ in degree 0 we add
a generator $\bar x_I$ in degree 1, and for each generator $y_I$
in degree 1 we add a generator $\bar y_I$ in degree 2. The
differential is as follows: $dy_I = x_I - \phi(x_I)$, as in
$U_{(T,\phi)}$, $d\bar x_I = x_I$, so we kill $\pi_0$, and $d\bar
y_I = y_I - \bar x_I - u_I$ where $u_I$ is a decomposable element
of degree 1 such that $d u_I = \phi(x_I)$. Such an element exists
indeed since $\phi(x_I)$ is a decomposable element in degree 0,
where all indecomposables are hit by a differential.


\begin{theorem}
\label{plus}

Let $\mathcal O$ be an admissible operad. Then the homotopical
nullification with respect to $\mathcal U$ is a functorial
plus-construction in the category of algebras over $\mathcal O$.
It enjoys the following properties:
\begin{itemize}
\item[]{(i)} $A^+ \simeq Cof(\coprod_{f\in Hom(\mathcal U,A)} \U
\longrightarrow A)$

\item[]{(ii)} $\pi_0(A^+)\cong\pi_0(A)/\mathcal P\pi_0(A)$

\item[]{(iii)} $H^Q_*(A)\cong H^Q_*(A^+)$
\end{itemize}
\end{theorem}

\begin{proof}
Consider the cofibre sequence
$$
\coprod_{f\in Hom(\mathcal U,A)} \mathcal U
\stackrel{ev}{\longrightarrow} A \longrightarrow B
$$
Clearly $A \rightarrow B$ is a $P_{\U}$-equivalence. So it remains
to show that $B$ is $\U$-local. Now, as the suspension of
$\mathcal U$ is contractible by the Hurewicz
Theorem~\ref{Hurewicz}, $B$ is $P_{\mathcal U}$-local if and only
if $[\U, B] = 0$. This happens exactly when $\mathcal P\pi_0(B) =
0$. Let us thus compute $\pi_0 B$. Consider actually the more
elementary cofiber $C_\alpha$ of a single map $\alpha: \mathcal
U_{(T, \phi)} \rightarrow A$. Such a map corresponds to an element
$a \in \mathcal P \pi_0 A$ together with a decomposition following
the pattern indicated by the tree $(T, \phi)$. Let us replace $A$
by a free algebra $\mathcal O (V)$ and construct now $C_\alpha$ as
the push-out of $\mathcal O (V) \leftarrow U_{(T, \phi)} \mono
C(U_{(T, \phi)})$. The models of these algebras we exhibited
earlier show that $C_\alpha = \mathcal O(V) \coprod \mathcal
O(\bar x_I, \bar y_I)$ with $d \bar x_I = a_I = \alpha(x_I)$ and
$d \bar y_I = b_I - \bar x_I - \alpha(u_I)$. Clearly $\pi_0
C_\alpha \cong \pi_0 A / <a>$. Likewise $\pi_0 B \cong \pi_0 A /
\mathcal P\pi_0(A)$.

Hence $\mathcal P\pi_0(B) = 0$, which shows that $B \simeq A^+$
and the other two properties are now direct consequences of the
first one.
\end{proof}

\medskip

\noindent {\bf Naturality.} We conclude this section with a
discussion of the naturality of the plus-construction with respect
to the operad. We denote by $\mathcal U'$ the universal acyclic
$\mathcal O'$-algebra as constructed above and $A^{+'} =
P_{\mathcal U'} A$ the associated plus-construction.

\begin{proposition}
\label{natural}

Let $\mathbf f:\mathcal O\longrightarrow \mathcal O'$ be a map of
operads, then there is a map of $\mathcal O$-algebras $f:\mathcal
U \longrightarrow \mathcal U'.$
\end{proposition}

\begin{proof}
The map $\mathbf f$ induces a map between the directed systems
$\{U_{r},\phi_r\}$ and $\{U_{r}',\mathbf f(\phi_r)\}$. Where
$\{U_{r},\phi_r\}$ is the directed system associated to a
$\mathcal O$-tree $(T,\phi)$ and $\{U_{r}',\mathbf f(\phi_r)\}$ is
the directed system associated to the $\mathcal O'$-tree
$(T,\mathbf f(\phi))$ where each vertex is of the form $\mathbf
f(o)$. There is a natural transformation between the directed
systems of $\mathcal O$-algebras, thus also a map between their
homotopy colimits.
\end{proof}

\begin{proposition}
\label{compatibility}

Let $\mathbf f:\mathcal O\longrightarrow \mathcal O'$ be a
quasi-isomorphism of operads, and suppose that either we work over
$\mathbb Q$, or the operads $\mathcal O$ and $\mathcal O'$ are
cofibrant. Then $f: \mathcal U \rightarrow \mathcal U'$ is a
quasi-isomorphism of $\mathcal O$-algebras.
\end{proposition}

\begin{proof}
The result follows from the fact that free algebras over the
operads $\mathcal O$ and $\mathcal O'$ and over the same
generators are quasi-isomorphic as $\mathcal O$-algebras.
\end{proof}

\medskip

As a consequence, when replacing an operad by a cofibrant one to
do homotopy, the choice of this cofibrant operad does not matter.

\begin{corollary}
\label{welldef}
Let $A$ be an $\mathcal O'$-algebra, and let $\mathbf f:\mathcal O
\longrightarrow \mathcal O'$ be a morphism of operads. Under the
same assumptions as in the preceding proposition, the map $A^+
\longrightarrow A^{+'}$ is a quasi-isomorphism of $\mathcal
O$-algebras.\hfill{$\square$}
\end{corollary}

\section{Fibrations and the plus-construction}

Let $\mathcal O$ be an admissible operad over a field $k$. This
section is devoted to the analysis of the behavior of the
plus-construction with fibrations. In particular we will be
interested in the fibre $AX$ of the map $X \rightarrow X^+$. As
one should expect it, $AX$ is the universal acyclic algebra over
$X$ in the sense that any map $A \rightarrow X$ from an acyclic
algebra $A$ factors through $AX$. The most efficient tool to deal
with such questions is the technique of fibrewise localization in
our model category of $\mO$-algebras. To our knowledge, such a
tool had not been developed up to now in any other context than
spaces, and we refer therefore to the separate paper \cite{DCS}
for the following claim:

\begin{theorem}
\label{fibrewise}
Let $F \rightarrow E \epi B$ be a fibration of $\mO$-algebras.
There exists then a commutative diagram
\[
\diagram F \rto \dto & E \rto \dto & B \dto\cr F^+ \rto & \bar E
\rto & B
\enddiagram
\]
where both lines are fibrations and the map $E \rightarrow \bar E$
is a $P_{\cal U}$-equivalence.\hfill{$\square$}
\end{theorem}

The main ingredient in the proof of this theorem is the fact that
the category of $\mO$-algebras satisfies the cube axiom. From the
above theorem we infer that the plus-construction sometimes
preserves fibrations.

\begin{theorem}
\label{preserve}
Let $F\longrightarrow E\longrightarrow B$ be a fibration of
$\mathcal O$-algebras. If the basis $B$ is local with respect to
the $\mathcal U$-nullification then we have a fibration
$$F^+\longrightarrow E^+\longrightarrow B.$$
\end{theorem}

\begin{proof}
By Theorem \ref{fibrewise} this is a direct consequence of the
fact that the total space $\bar E$ sits in a fibration where both
the fibre and the base space are $\mathcal U$-local and hence is
also $\mathcal U$-local.
\end{proof}
\\
\\
\noindent {\bf The fibre of the plus-construction.} Another
consequence of the fibrewise plus-construction is that the fibre
$AX$ is acyclic.


\begin{proposition}
\label{acyclic}
The fiber $AX$ of the plus-construction $X\longrightarrow X^+$ is
$H^Q_*$-acyclic for any $\mO$-algebra $X$.
\end{proposition}

\begin{proof}
Consider the fibration $AX \rightarrow X \rightarrow X^{+}$. The
plus-construction preserves this fibration by the above theorem,
i.e. $(AX)^{+}$ is the fibre of the identity on $X^{+}$. Thus
$(AX)^{+} \simeq *$, which means that $AX$ has the same Quillen
homology as a wedge of copies of $\mathcal U$.
\end{proof}

\medskip

\noindent {\bf Cellularization.} We can go a little further in the
analysis of the fibre $AX$. Our next result says precisely that
the map $AX \rightarrow X$ is a $CW_{\mathcal U}$-equivalence,
where $CW_{\mathcal U}$ is Farjoun's cellularization functor
(\cite[Chapter~2]{Dro}). We do not know whether $AX$ is actually
the $\mathcal U$-cellularization of~$X$.


\begin{proposition}
\label{cellular}

We have $map(\mathcal U, AX) \simeq map(\mathcal U, X)$ for any
$\mO$-algebra $X$.
\end{proposition}

\begin{proof}
We know $AX$ is acyclic by Proposition~\ref{acyclic}. Apply now
$map(\mathcal U, -)$ to the fibration $AX \rightarrow X
\rightarrow X^+$ so as to get a fibration of simplicial sets
$$
map(\mathcal U, AX) \rightarrow map(\mathcal U, X) \rightarrow
map(\mathcal U, X^+)
$$
By definition $X^+$ is $\mathcal U$-local, so that $map(\mathcal
U, X^+)$ is trivial. Therefore $map(\mathcal U, AX) \simeq
map(\mathcal U, X)$.
\end{proof}

\medskip

On the level of components, this implies we have an isomorphism
$[\mathcal U, AX] \cong [\mathcal U, X]$, which means that any
element in the $\mO$-perfect ideal $\mathcal P \pi_0 X$ together
with a given decomposition can be lifted in a unique way to such
an element in $\pi_0 AX$.

\begin{proposition}
\label{cofibration}
The fibration $AX \rightarrow X \rightarrow X^+$ is also a
cofibration.
\end{proposition}

\begin{proof}
By definition $X^+$ is the homotopy cofibre of a map $\coprod
\mathcal U \rightarrow X$. By the above proposition this map
admits a unique lift to $AX$. By considering the composite
$\coprod \mathcal U \rightarrow AX \rightarrow X$, we get a
cofibration
$$
Cof(\coprod \mathcal U \rightarrow AX) \rightarrow Cof(\coprod
\mathcal U \rightarrow X) \rightarrow Cof(AX \rightarrow X)
$$
The first cofibre is $(AX)^+$, which is contractible, and the
second is $X^+$. The third is thus $X^+$ as well.
\end{proof}

\medskip

\noindent {\bf Preservation of square zero extensions.} Let us
finally study the effect of the plus-construction on a square zero
extension, as introduced at the end of Section~1. In the case of
Lie or Leibniz algebras this notion coincides of course with the
classical one of central extension, as exposed e.g. in~\cite{KL}.
Following~\cite{QU}, \cite [chapter 5]{GH}, such a square zero
extension is classified by an element in the first Quillen
cohomology group $H^1_Q(B; Z) \cong [B, K(Z, 1)]$. Here $K(Z, 1)$
is the delooping of $Z$, given as $\mO$-algebra by the chain
complex $Z$ concentrated in degree~1. As for group extensions, the
fibre of the classifying map $B \rightarrow K(Z, 1)$ (the
$k$-invariant of the extension) is precisely $E$.

\begin{proposition}
\label{central}
Let $Z \hookrightarrow E \rightarrow B$ be a square zero extension
of discrete $\mO$-algebras. Then $Z \rightarrow E^+ \rightarrow
B^+$ is a fibration.
\end{proposition}

\begin{proof}
Let us consider the $k$-invariant and the associated fibration $E
\rightarrow B \rightarrow K(Z, 1)$. The base is 0-connected, thus
$P_{\mathcal U}$-local. Theorem~\ref{preserve} tells us that $E^+
\rightarrow B^+ \rightarrow K(Z, 1)$ is also a fibration.
Therefore so is $Z \rightarrow E^+ \rightarrow B^+$.
\end{proof}

\begin{section}{Applications to algebras of matrices}

{\bf Recollections on algebras of matrices.} Let $k$ be a field
and $R$ be an associative $k$-algebra. Consider $gl(R)$ the union
of the $gl_n(R)$'s. This is a $\mathcal Lie$-algebra and also a
$\mathcal Leib$-algebra for the classical bracket of matrices. The
trace $tr: gl(R)\longrightarrow R/[R,R]$ is a morphism of
$\mathcal Lie$ and $\mathcal Leib$-algebras, whose kernel is by
definition the algebra $sl(R)$.
\\
We define the Steinberg algebra $st(R)$ for the two operads
$\mathcal Lie$ and $\mathcal Leib$ by taking the free algebras
over the generators $u_{i,j}(r)$, $r\in R$ and $1\leq i\not = j$
with the relations
\\
a) $u_{i,j}(m.r+n.s)=m.u_{i,j}(r)+n.u_{i,j}(s)$ for $r,s\in R$ and
$m,n\in \mathbb Z$.
\\
b) $[u_{i,j}(r),u_{k,l}(s)]=0$ if $i \not =l$ and $j\not =k$.
\\
c) $[u_{i,j}(r),u_{k,l}(s)]=u_{i,l}(rs)$ if $i \not =l$ and $j=k$.
\\
We have the following extension of algebras (for both operads):
$$
Z(st(R))\longrightarrow st(R) \longrightarrow sl(R)
$$
where $Z(st(R))$ is the kernel of the canonical map between
$st(R)$ and $sl(R)$. Following the work of C.~Kassel and
J.L.~Loday this is a universal square zero extension
\cite[Proposition~1.8]{KL}.
\\
Now we can consider all these algebras as algebras over cofibrant
replacements $\mathcal L_{\infty}$ and $\mathcal Leib_{\infty}$ of
the operads $\mathcal Lie$ and $\mathcal Leib$.
\\
\\
{\bf Homology theories.} In the category of $\mathcal
L_{\infty}$-algebras we define cyclic homology up to homotopy
$HC^{\infty}$:
$$
HC_{*}^{\infty}(R)=\pi_*(gl(R)^+).
$$
Likewise in the category of $\mathcal Leib_{\infty}$-algebras we
define Hochschild homology up to homotopy:
$$
HH_{*}^{\infty}(R)=\pi_*(gl(R)^+).
$$
By Corollary~\ref{welldef}, we notice that these definitions do
not depend on the choice of the cofibrant replacement of the
operads $\mathcal Lie$ or $\mathcal Leib$. These theories define
two functors from the category of associative algebras to the
categories of $\mathcal Lie$ and $\mathcal Leib$ graded algebras.
We recall that the homotopy of a $\mathcal L_{\infty}$-algebra
(resp. a $\mathcal Leib_{\infty}$-algebra) is a graded $\mathcal
Lie$-algebra (resp. a $\mathcal Leib$-algebra).
\\
When we consider these two theories over $\mathbb Q$, we have
quasi-isomorphisms $\mathcal Leib_{\infty}\rightarrow \mathcal
Leib$ and $\mathcal Lie_{\infty}\rightarrow \mathcal Lie$. Then by
the results of M. Livernet \cite{ML} our theories coincide with
the classical cyclic and Hochschild homologies:
$$HC^{\infty}(R)\cong HC(R),$$
$$HH^{\infty}(R)\cong HH(R).$$
\\
We do not know if these isomorphisms remain valid over $\mathbb
Z$. However, using the properties of our construction, we are able
to compute the first four groups of $HC^{\infty}$ and
$HH^{\infty}$. These results form perfect analogues of the
classical computations in algebraic $K$-theory, see for example
\cite[Theorem~4.2.10]{Ros}, and \cite[Theorem~3.14]{Arl} for a
topological approach.
\\
\\
{\bf Abelianization.} In order to compute $HH_0^{\infty}(R)$ and
$HC_0^{\infty}(R)$ we use the following fibration given by the
trace:
$$
sl(R)\longrightarrow gl(R)\longrightarrow R/[R,R].
$$

\begin{proposition}
\label{K0}
Let $R$ be an associative $k$-algebra. Then $HH_0^{\infty}(R)$ and
$HC_0^{\infty}(R)$ are both isomorphic to $R/[R,R]$.
\end{proposition}

\begin{proof}
By Theorem~\ref{preserve} we get a fibration
$$
sl(R)^+\longrightarrow gl(R)^+\longrightarrow R/[R,R].
$$
The commutator subgroup of $gl(R)$ as well as $sl(R)$ (i.e. the
perfect radical in either the category of Lie or Leibniz algebras)
is $sl(R)$. Therefore so is the perfect radical in $\mathcal
L_\infty$ and $\mathcal Leib_\infty$ (this is the case for any
discrete algebra). Hence $\pi_0 gl(R)^+ \cong R/[R, R]$.
\end{proof}
\\
\\
{\bf The center of the Steinberg algebra.} In order to
compute $HH_1^{\infty}(R)$ and $HC_1^{\infty}(R)$, we use the
Steinberg Lie, respectively the Steinberg Leibniz, algebra $st(R)$
and the following square zero extension:
$$
Z(st(R))\longrightarrow st(R) \longrightarrow sl(R).
$$
This is the universal central extension of the perfect algebra
$sl(R)$. In particular $st(R)$ is superperfect, meaning that
$H_1^Q(st(R)) = 0$.

\begin{proposition}
\label{K1}
Let $R$ be an associative $k$-algebra. Then the first homology
group $HC_1^{\infty}(R)$ (respectively $HH_1^{\infty}(R)$) is
isomorphic to $Z(st(R)) \cong H_1^Q(sl(R))$ (respectively to the
center of the Steinberg Leibniz algebra).
\end{proposition}

\begin{proof}
We have to compute $\pi_1 gl(R)^+$. From Theorem~\ref{preserve} we
infer that $sl(R)^+ \rightarrow gl(R)^+ \rightarrow R/[R, R]$ is a
fibration. Hence the preceding proposition tells us that $sl(R)$
is the $0$-connected cover of $gl(R)^+$, so that we only need to
compute $\pi_1 sl(R)^+$. By the Hurewicz Theorem, this is
isomorphic to $H_1^Q(sl(R))$.

Moreover Proposition~\ref{central} shows that $Z(st(R))
\rightarrow st(R)^+ \rightarrow sl(R)^+$ is a fibration. Both
$sl(R)$ and $st(R)$ are perfect algebras, so their
plus-constructions are $0$-connected. Actually $st(R)^+$ is even
$1$-connected since $H_1^Q(st(R)) = 0$. The homotopy long exact
sequence allows now to conclude that $\pi_1 sl(R)^+ \cong
Z(st(R))$.
\end{proof}

\begin{proposition}
\label{K2}
Let $R$ be an associative $k$-algebra. Then for $2 \leq i \leq 3$,
$HC_i^{\infty}(R)$ (respectively $HH_i^{\infty}(R)$) is isomorphic
to $H_i^Q(st(R))$, the Quillen homology of the Steinberg Lie
algebra in the category of $\mathcal L_{\infty}$-algebras
(respectively to $H_2^Q(st(R))$, the Quillen homology of the
Steinberg Leibniz algebra in the category of $\mathcal
Leib_{\infty}$-algebras).
\end{proposition}

\begin{proof}
This is the same proof as the first part of the computation of
$HH_2^\infty$ and $HC_2^\infty$. The computation of the third
groups follows from the Hurewicz Theorem \ref{Hurewicz}.
\end{proof}

\medskip

As explained in the first section there is an isomorphism over
$\mathbb Q$ between the Quillen homology $H^Q_*$ and
$H^{Lie}_{*+1}$, respectively $H^{Leib}_{*+1}$. Together with the
fact that the theories up to homotopy coincide with their
classical analogues over $\mathbb Q$, the three computations we
made above yield the following isomorphisms.

\begin{corollary}
\label{classical}
 Let $R$ be an associative algebra over $\mathbb Q$ then
\begin{itemize}
\item[(1)] $HC_0(R)=HH_0(R)=R/[R,R]$, \item[(2)]
$HC_1(R)=H_2^{Lie}(sl(R))$, $HH_1(R)=H_2^{Leib}(sl(R))$,
\item[(3)] $HC_2(R)=H_3^{Lie}(st(R))$,
$HH_2(R)=H_3^{Leib}(st(R))$,\item[(4)] $HC_3(R)=H_4^{Lie}(st(R))$,
$HH_3(R)=H_4^{Leib}(st(R))$. \hfill{$\square$}
\end{itemize}
\end{corollary}

\noindent{\bf Morita invariance.} These theories are obviously
Morita invariant since $gl(gl(R))$ is isomorphic to $gl(R)$. Hence
we have $HC_*^{\infty}(gl(R))\cong HC_*^{\infty}(R)$ and
$HH_*^{\infty}(gl(R))\cong HH_*^{\infty}(R)$.

\medskip

\noindent{\bf Products.} Let $R$ and $S$ be two associative
$k$-algebras, and form the product in the category of associative
algebras $R\times S$. We want to compute $HC_*^{\infty}(R\times
S)$ and $HH_*^{\infty}(R\times S)$. Observe that $gl(R\times S)$
is isomorphic as a Lie-algebra to the product $gl(R)\times gl(S)$.
As nullifications preserve products (this is a consequence of the
fiberwise localization \cite{DCS}) one has:

\begin{proposition}
Let $R$ and $S$ be two associative $k$-algebras. Then:
\begin{itemize}
\item[(i)]{} $HC_*^{\infty}(R\times S)\cong HC_*^{\infty}(R)\oplus
HC_*^{\infty}(S)$
\item[(ii)]{} $HH_*^{\infty}(R\times S)\cong HH_*^{\infty}(R)\oplus HH_*^{\infty}(S)$ \hfill{$\square$}
\end{itemize}
\end{proposition}

\end{section}

\bibliographystyle{amsplain}

\vskip 0.5 cm

\setlength{\baselineskip}{0.6cm}

\bigskip\noindent
David Chataur:

 \noindent Centre de Recerca Matematica, E--08193 Bellaterra, Spain,
 email: {\tt dchataur@crm.es}

\bigskip\noindent
Jos\'e L. Rodr\'{\i}guez:

\noindent Departamento de Geometr\'\i a, Topolog\'\i a y Qu\'\i mica
Org\'anica, Universidad de Almer\'\i a, E--04120 Almer\'\i a, Spain,
e-mail: {\tt jlrodri@ual.es}

\bigskip\noindent
J\'er\^ome Scherer:

\noindent Departament de matem\`atiques, Universitat Aut\'onoma de
Barcelona, E--08193 Bellaterra, Spain, e-mail: {\tt
jscherer@mat.uab.es}

\end{document}